\documentclass[11pt,leqno]{article}
\usepackage[centertags]{amsmath}
\usepackage{amsfonts}
\usepackage{amsmath}
\usepackage{amssymb}
\usepackage{latexsym}
\usepackage{amsthm}
\usepackage{mathrsfs}
\usepackage{newlfont}
\usepackage[latin1]{inputenc}

\usepackage{epsfig}
\usepackage{pst-all}
\usepackage{graphicx}
\usepackage{pstricks}

\newtheorem{teor}{Theorem}[section]

\newtheorem{nota}{Remark}[section]

\numberwithin{equation}{section} 

\title{A note on the invariant distribution of a quasi-birth-and-death process\thanks{The work of the author is partially supported by D.G.E.S, ref. BFM2006-13000-C03-01, Junta de Andaluc\'{i}a, grants FQM-229, FQM-481, P06-FQM-01738 and Subprograma de estancias de movilidad posdoctoral en el extranjero, MICINN, ref. -2008-0207.\newline \textbf{AMS Subject
Classifications.} 60J10, 42C05. \newline \textbf{Key words.} Quasi-birth-and-death processes, matrix-valued orthogonal polynomials, Markov chains, block tridiagonal transition probability matrix}}

\author{Manuel D. de la Iglesia \\ \footnotesize Courant Institute of Mathematical Sciences. New York University \\ \footnotesize\ 251 Mercer Street, New York, NY 10012, U.S.A. mdi29@cims.nyu.edu}

\date{}

\begin{document}

\maketitle

\begin{abstract}
The aim of this paper is to give an explicit formula of the invariant distribution of a quasi-birth-and-death process in terms of the block entries of the transition probability matrix using a matrix-valued orthogonal polynomials approach. We will show that the invariant distribution can be computed using the squared norms of the corresponding matrix-valued orthogonal polynomials, no matter if they are or not diagonal matrices. We will give an example where the squared norms are not diagonal matrices, but nevertheless we can compute its invariant distribution.
\end{abstract}

\section{Introduction}\label{INTRO}

The connection between random walks/birth-and-death processes and orthogonal polynomials is very well known. In both cases the state space is the set of nonnegative integers, i.e. $\mathcal{S}=\{0,1,\ldots\}$, while the parameter set is $\mathcal{T}=\{0,1,\ldots\}$ for random walks (i.e. discrete time) and $\mathcal{T}=[0,\infty)$ for birth-and-death processes (i.e. continuous time). In a series of papers, Karlin and McGregor (see \cite{KMc1, KMc2, KMc3}) found an appropriate tool to study these processes connecting the tridiagonal one-step transition probability matrix $P$ (or the infinitesimal generator $\mathcal{A}$ in the continuous time case) with a measure supported in the real line. In particular, they obtained an integral representation for the $n$-step transition probability matrix $P^n$ (or the transition probability matrix $P(t)$ in the continuous case) in terms of the spectral measure and the corresponding orthogonal polynomials, as well as the explicit expression of the invariant distribution, i.e. a row vector $\mbox{\boldmath$\pi$}$ with nonnegative components such that $\mbox{\boldmath$\pi$}P=\mbox{\boldmath$\pi$}$. 

Quasi-birth-and-death processes are a natural extension where now the state space is two dimensional, i.e. $\mathcal{S}=\{(i,j) : i\geq0, 1\leq j\leq N\}$. The first component of the pair is usually called the \emph{level} and the second one the \emph{phase}. Now the one-step transition probability matrix $P$ is block tridiagonal. The link with matrix-valued orthogonal polynomials was initially raised independently by \cite{DRSZ} and \cite{G2} in discrete time (for continuous time see \cite{DR}), although the first example may be traced back to the last section of \cite{KMc3}, where the authors deal with the case of a random walk with state space the set of all integers, replacing the spectral measure by a $2 \times 2$ non-negative matrix. In \cite{DRSZ} and \cite{G2} one can find an integral representation for the $n$-step transition probability matrix $P^n$, along with other probabilistic useful results. For a much more detailed presentation of this field, as well as its connections with queueing problems in network theory and the general field of communication systems the reader should consult \cite{LR1, N, O}.

Nevertheless computing the invariant distribution of a quasi-birth-and-death process, i.e. a row vector $\mbox{\boldmath$\pi$}$ with nonnegative components such that $\mbox{\boldmath$\pi$}P=\mbox{\boldmath$\pi$}$, is a much harder problem, compared with the 1-dimensional situation. It is possible to derive nicer looking expressions for the invariant distribution for some special cases, like level-independent quasi-birth-and-death processes (see \cite{LPT}). Also in some papers (see \cite{G3, G4, G5, GdI2}) the invariant distribution was given in terms of the entries of the \emph{diagonal} norms of the corresponding matrix-valued orthogonal polynomials, giving a natural candidate. Here we will show that it is possible to derive an explicit expression of the invariant distribution in terms of the block entries of the transition probability matrix or equivalently, in terms of the squared norms of the corresponding matrix-valued orthogonal polynomials, no matter if they are or not diagonal. In fact, we will give an example of a quasi-birth-and-death process where the squared norms of the corresponding matrix-valued orthogonal polynomials are \emph{not} diagonal, but nevertheless we can compute its invariant distribution.

In Section \ref{QBDP} we will define the matrix version of the so-called potential coefficients and connect them with the squared norms of matrix-valued orthogonal polynomials. In Section \ref{IDs} we will give the explicit expression of an invariant distribution in terms of these matrix-valued potential coefficients and finally, in Section \ref{EX}, we will give the example mentioned in the paragraph above.

\section{Matrix-valued potential coefficients}\label{QBDP}

The content in this section is already known in the literature, but we will point out some properties that these matrix-valued potential coefficients have in order to prove the main result in the next section. For simplicity, we will focus on the case of discrete time quasi-birth-and-death processes. The continuous time case only suffers cosmetic changes.

Consider a nonhomogeneous discrete time quasi-birth-and-death process, i.e. a discrete time Markov process on the countable two dimensional state space $\mathcal{S}=\{(i,j) : i\geq0, 1\leq j\leq N\}$ with block tridiagonal transition probability matrix
\begin{equation}\label{PQdb}
P=\begin{pmatrix}
  B_0 & A_0   \\
  C_1 & B_1 & A_1 &  \\
  & C_2 & B_2 & A_2  \\
   &  & \ddots & \ddots & \ddots
\end{pmatrix}.
\end{equation}

The first component of the pair is usually called the level and the second one the phase. $P$ is \emph{stochastic}, i.e. all entries are nonnegative and all rows sum up to one, that is $(B_0+A_0)\textbf{e}_N=\textbf{e}_N$ and $(C_n+B_n+A_n)\textbf{e}_N=\textbf{e}_N$, $n\geq1$, where $\textbf{e}_N$ is the column vector of 1's of dimension $N$. We will assume that our process is irreducible and that the coefficients $A_n$ and $C_n$ are nonsingular matrices.  Clearly in the case when the number of phases $N$ is one we are back to the case of an ordinary birth-and-death process.

In order to link these processes with matrix-valued orthogonal polynomials the matrix $P$ needs to be transformed into a symmetric matrix. This is possible if there exists a nonsingular block diagonal matrix
$$
R=\begin{pmatrix}
  R_0 &    \\
   & R_1 &    \\
     &  &\ddots
\end{pmatrix}
$$
such that $RPR^{-1}$ is symmetric (see for instance Theorem 2.1 of \cite{DRSZ}). The matrices $R_n$ are subject to the following restrictions
\begin{equation}\label{RnR}
R_nB_nR_n^{-1}=(R_nB_nR_n^{-1})^T,\quad\mbox{and}\quad R_nA_nR_{n+1}^{-1}=(R_{n+1}C_{n+1}R_n^{-1})^T,\quad n\geq0,
\end{equation}
where $M^T$ denotes the transpose of the matrix $M$.

Let us call
\begin{equation}\label{Pin}
\Pi_n=R_n^TR_n,\quad n\geq0,
\end{equation}
which are obviously symmetric. From \eqref{RnR} we have that $\Pi_n$, $n\geq0$, satisfy
\begin{equation}\label{Con1}
\Pi_nB_n=B_n^T\Pi_n,\quad n\geq0,
\end{equation}
and
\begin{equation}\label{Con2}
\Pi_nA_n=C_{n+1}^T\Pi_{n+1},\quad n\geq0.
\end{equation}

Condition \eqref{Con2} gives an explicit formula for $\Pi_n$, $n\geq1$, given $A_n$, $C_n$ and $\Pi_0$:
\begin{equation}\label{MVPC}
\Pi_n=(C_1^TC_2^T\cdots C_{n}^{T})^{-1}\Pi_{0}(A_0A_1\cdots A_{n-1}),\quad n\geq1,
\end{equation}
This formula is similar to the scalar potential coefficients for birth-and-death processes (see \cite{KMc2}), so we will call them \emph{matrix-valued potential coefficients}.

Under these assumptions, there always exists a weight matrix $W$ such that the matrix-valued polynomials defined by the three-term recurrence relation
\begin{equation}\label{ttrrst}
xQ_n(x)=A_{n}Q_{n+1}(x)+B_nQ_n(x)+C_{n}Q_{n-1}(x),\quad n\geq0,
\end{equation}
where $Q_{-1}(x)=0$ and $Q_0(x)=I$, are orthogonal, i.e.
$$
\int Q_n(x)W(x)Q_m^T(x)dx=\|Q_n\|^2_W\delta_{nm}.
$$
In particular, as it is remarked in page 121 of \cite{DRSZ}, $\Pi_0$ can be given in terms of the 0-th moment $S_0=\int W(x)dx$. For an orthonormal family $(\widetilde{Q}_n)_n$ we have that
$$
I=\int \widetilde{Q}_0(x)W(x)\widetilde{Q}_0^T(x)dx=R_0S_0R_0^{T}
$$
Therefore
\begin{equation*}\label{Pn0}
\Pi_0=S_0^{-1}=(\|Q_0\|_W^2)^{-1}.
\end{equation*}
That means that for every weight matrix $W$ associated with the Jacobi matrix \eqref{PQdb}, the matrix-valued potential coefficients are defined recursively using \eqref{MVPC} in a unique way.

\begin{nota}
 Observe that the sequence $(R_n)_n$ is not unique in the sense that we can always consider any other sequence $(U_nR_n)_n$, where $(U_n)_n$ is any sequence of orthogonal matrices and we will get another family of orthonormal polynomials. Nevertheless, the matrix-valued potential coefficients $\Pi_n$ are always unique (see \eqref{Pin}).
\end{nota}

Finally we will show that the matrix-valued potential coefficients can be given in terms of the squared norms of the matrix-valued orthogonal polynomials. From the three-term recurrence relation \eqref{ttrrst} we can multiply on the right by $WQ_n$, $WQ_{n+1}$ and $WQ_{n-1}$, respectively, and integrate in such a way that the sequence of squared norms $\|Q_n\|_W^2$ satisfies the following relations
$$
B_n\|Q_n\|_W^2=\|Q_n\|_W^2B_n^T,\quad \|Q_n\|_W^2C_{n+1}^T=A_n\|Q_{n+1}\|_W^2,\quad n\geq0.
$$
But these are just the relations \eqref{Con1} and \eqref{Con2}. Since $\Pi_0=(\|Q_0\|_W^2)^{-1}$, we have that
\begin{equation}\label{Pn0N}
\Pi_n=(\|Q_n\|_W^2)^{-1},\quad n\geq0.
\end{equation}

In particular we have that the matrix-valued potential coefficients $\Pi_n$ are positive semi-definite matrices.

\section{The invariant distribution}\label{IDs}
With all the information from the last section we are able to give a very natural candidate for the invariant distribution of the process $P$ using the matrix-valued potential coefficients $\Pi_n$.

\begin{teor}\label{Teo}
Let $P$ be the transition probability matrix given by \eqref{PQdb}. Define the sequence of matrices $\Pi_n$, $n\geq1$, as in \eqref{MVPC} with $\Pi_0=(\int W(x)dx)^{-1}$ or, equivalently, as in \eqref{Pn0N}. Consider the following row vector
\begin{equation}\label{ID}
\mbox{\boldmath$\pi$}=
((\Pi_0\textbf{e}_N)^T;(\Pi_1\textbf{e}_N)^T;(\Pi_2\textbf{e}_N)^T;\cdots),
\end{equation}
where $\textbf{e}_N$ is the column vector of 1's of dimension $N$.  Then $\mbox{\boldmath$\pi$}$ is an invariant distribution for the process $P$, i.e. all components of $\mbox{\boldmath$\pi$}$ are nonnegative and
\begin{equation}\label{IDprop}
\mbox{\boldmath$\pi$}P=\mbox{\boldmath$\pi$}.
\end{equation}
\end{teor}
\proof
All components of $\mbox{\boldmath$\pi$}$ are nonnegative since $\Pi_n$ are positive semi-definite matrices (see \eqref{Pn0N}).

To prove \eqref{IDprop}, we have to check that
$$
(\Pi_0\textbf{e}_N)^TB_0+(\Pi_1\textbf{e}_N)^TC_1=(\Pi_0\textbf{e}_N)^T
$$
and
$$
(\Pi_{n-1}\textbf{e}_N)^TA_{n-1}+(\Pi_{n}\textbf{e}_N)^TB_n+(\Pi_{n+1}\textbf{e}_N)^TC_{n+1}=(\Pi_n\textbf{e}_N)^T,\quad n\geq1.
$$
The first equality holds using properties \eqref{Con1}, \eqref{Con2} and the fact that $P$ is stochastic, since
$$
\textbf{e}_N^T(\Pi_0B_0+\Pi_1C_1)=\textbf{e}_N^T(B_0^T\Pi_0+A_0^T\Pi_0)=[(B_0+A_0)\textbf{e}_N]^T\Pi_0=(\Pi_0\textbf{e}_N)^T,
$$
Similarly, for $n\geq1$
$$
\textbf{e}_N^T(\Pi_{n-1}A_{n-1}+\Pi_{n}B_n+\Pi_{n+1}C_{n+1})=[(A_n+B_n+C_n)\textbf{e}_N]^T\Pi_n=(\Pi_n\textbf{e}_N)^T.
$$
\hfill \endproof

\begin{nota}
The same result holds in the continuous time case. In this case, the transition probability matrix $P$ is replaced by an infinitesimal generator $\mathcal{A}$, block tridiagonal as in \eqref{PQdb}, but with $(B_0+A_0)\textbf{e}_N=0$ and $(C_n+B_n+A_n)\textbf{e}_N=0$, $n\geq1$ (see \cite{DR} for more details). Therefore the row vector $\mbox{\boldmath$\pi$}$ defined by \eqref{ID} is an invariant distribution for the process, i.e. all components are nonnegative and
$$
\mbox{\boldmath$\pi$}\mathcal{A}=0.
$$
\end{nota}

\bigskip

Let us make a few comments about the unicity of this invariant distribution. If the process is recurrent then there exists a unique invariant distribution given by \eqref{ID} (see Theorem 5.4 of \cite{S}). But if the process is transient there are no results about unicity. In fact, the case of a random walk on the integers
treated in \cite{KMc3} gives (for $p\neq q$) an example of a transient process where the invariant distribution is not unique.

In Section \ref{QBDP} we saw that the matrix-valued potential coefficients $\Pi_n$ are defined in terms of the corresponding weight matrix $W$ and the block entries of $P$. These quantities are constructed in a unique way, so the question of the unicity of the invariant distribution is the same as the question of the unicity of the weight matrix $W$ corresponding to the equivalent class  of block tridiagonal Jacobi matrices $P$ by unitary transformations. So if we can guarantee that the corresponding weight matrix is unique, then the invariant distribution will be also unique. Nevertheless there are not general results to determine if the weight matrix associated with a block tridiagonal Jacobi matrix is unique.

%

\bigskip

Observe that Theorem \ref{Teo} gives a way to compute an invariant distribution even if the squared norms of the corresponding matrix-valued orthogonal polynomials are not diagonal. Until now, all the examples introduced in the literature had the special lucky property that these norms were diagonal matrices. In the following section we will give an example of a quasi-birth-and-death process where the squared norms of the corresponding matrix-valued orthogonal polynomials are \emph{not} diagonal, but nevertheless, using \eqref{ID}, we can compute its invariant distribution. This example, as far as the author knows, may be the first quasi-birth-and-death process with this property in the literature.

\section{The example}\label{EX}

Consider the following $2\times2$ weight matrix
\begin{equation*}\label{W}
 W(x)=x^{\alpha}(1-x)^{\beta}\begin{pmatrix}
                                   kx^2+\beta-k+1 & (\beta-k+1)(1-x)\\
                                   (\beta-k+1)(1-x) & (\beta-k+1)(1-x)^2
                                   \end{pmatrix},\quad x\in[0,1],
\end{equation*}
where $\alpha,\beta>-1$ and $0<k<\beta+1$. This weight matrix was introduced in \cite{DdI2} in a different context, that of searching second-order differential operators having several families of matrix-valued orthogonal polynomials as eigenfunctions. Let us call $(\widehat{P}_n)_n$ the unique monic family of matrix-valued orthogonal polynomials associated with $W$.

Let $\Delta_0$ be the nonsingular matrix
\begin{equation*}\label{TT}
\Delta_0=\begin{pmatrix}
                                   1 & -\frac{\alpha+\beta-k+3}{\alpha+2\beta-2k+4}\\
                                   1 & -\frac{\alpha+2\beta-2k+4}{\beta-k+1}
                                   \end{pmatrix},
\end{equation*}
and let us consider the equivalent weight matrix
\begin{equation*}\label{Ww}
\widetilde{W}=\Delta_0W\Delta_0^T.
\end{equation*}
Now we consider a special family of matrix-valued orthogonal polynomials with respect to $\widetilde{W}$ of the following form
\begin{equation*}\label{Qn}
Q_n(x)=\Delta_n\widehat{P}_n(x)\Delta_0^{-1},
\end{equation*}
where
\begin{equation*}\label{Deln}
\Delta_n=\begin{pmatrix}
                             {\footnotesize \mbox{$-\dfrac{a_n(n^2+n(\alpha+\beta+3)+k(\alpha+2\beta-2k+4))}{k(\alpha+\beta-k+n+3)}$}} &a_n\\
                                  b_n  &{\footnotesize \mbox{$-\dfrac{b_n(n^2+n(\alpha+\beta+3)+k(\alpha+2\beta-2k+4))}{(k+n)(\beta-k+1)}$}}
                                   \end{pmatrix},
\end{equation*}
with
\begin{equation*}
a_n= {\Large \mbox{$-\frac{(\alpha+\beta+n+3)_n(\alpha+\beta-k+n+3)(n(\alpha+\beta+n+2)+k(\alpha+1))}{(\beta+2)_n(n(\alpha+\beta)^2+n(2n+5)(\alpha+\beta)+n(n+2)(n+3)+k(2\alpha\beta+2\beta-k(n+2)+\alpha^2+5\alpha-2k\alpha-n^2+4))}$}},
\end{equation*}
$$
b_n=\frac{(\alpha+\beta+n+3)_n(k+n)(n(\alpha+\beta+n+2)+k(\alpha+1))}{(\beta+2)_n((n^2+2nk)(\alpha+\beta)+n(n+5)k+k^2(\alpha-n+1)+n^2(n+3))},
$$
and $(z)_n$ will denote the Pochhammer symbol defined by
$(z)_{n}=z(z+1)\cdots(z+n-1)$ for $n>0$, $(z)_{0}=1$. Observe that $Q_0=I$. The choice of the leading coefficient $\Delta_n\Delta_0^{-1}$ of the family of polynomials $(Q_n)_n$ is motivated by the remarkable fact that
\begin{equation}\label{qn1}
    Q_n(1)\textbf{e}_2=\textbf{e}_2,
\end{equation}
where $\textbf{e}_2$ is the $2$-dimensional column vector with all
entries equal to 1. In other words, the sum of the elements in each
row of $Q_n(1)$ gives the value 1.

The family of matrix-valued orthogonal polynomials introduced above satisfies a
three-term recursion relation
\begin{equation}\label{ttrrst2}
xQ_n(x)=A_{n}Q_{n+1}(x)+B_nQ_n(x)+C_{n}Q_{n-1}(x),\quad n\geq0,
\end{equation}
where $Q_{-1}(x)=0$ and $Q_0(x)=I$ and $A_n, B_n, n\geq0$, $C_n, n\geq1$ are $2\times2$ full matrices with nonnegative entries. The corresponding Jacobi matrix
\begin{equation}\label{PP}
P=\begin{pmatrix}
  B_0 & A_0   \\
  C_1 & B_1 & A_1   \\
   & C_2 & B_2 & A_2   \\
   &  & \ddots & \ddots & \ddots
\end{pmatrix}
\end{equation}
is stochastic. This is a consequence of choosing the family $(Q_n)_n$ with the property \eqref{qn1}. Indeed, applying
both sides of the identity \eqref{ttrrst2} to the vector $\textbf{e}_2$, setting $x=1$ and using \eqref{qn1} we obtain that $(B_0+A_0)\textbf{e}_2=\textbf{e}_2$ and $(C_n+B_n+A_n)\textbf{e}_2=\textbf{e}_2$, $n\geq1$, i.e. the sum of the entries in each row of $P$ equals one. Therefore, it gives a quasi-birth-and-death process with two phases ($N=2$) and depending on three parameters, $\alpha, \beta$ and $k$.

The state space and the corresponding one-step transitions look as follows (see \cite{G2, G3, GdI2} for an explanation of the labeling and ordering of the states)

\bigskip

\bigskip

\bigskip

$$\begin{psmatrix}[rowsep=2.5cm,colsep=3cm]
  \cnode{.55}{0}& \cnode{.55}{2} & \cnode{.55}{4} & \cnode{.55}{6} & \pnode{8} \\
  \cnode{.55}{1} & \cnode{.55}{3} & \cnode{.55}{5} & \cnode{.55}{7} & \pnode{9} \\
\psset{nodesep=3pt,arcangle=15,labelsep=2ex,linewidth=0.3mm,arrows=->,arrowsize=1mm
3} \nccurve[angleA=130,angleB=170,ncurv=4]{0}{0}
\nccurve[angleA=190,angleB=230,ncurv=4]{1}{1}
\nccurve[angleA=70,angleB=110,ncurv=4]{2}{2}
\nccurve[angleA=70,angleB=110,ncurv=4]{4}{4}
\nccurve[angleA=70,angleB=110,ncurv=4]{6}{6}
\nccurve[angleA=250,angleB=290,ncurv=4]{3}{3}
\nccurve[angleA=250,angleB=290,ncurv=4]{5}{5}
\nccurve[angleA=250,angleB=290,ncurv=4]{7}{7} \ncarc{0}{2}
\ncarc{2}{0} \ncarc{2}{4} \ncarc{4}{2} \ncarc{4}{6} \ncarc{6}{4}
\ncarc{6}{8} \ncarc{8}{6} \ncarc{0}{1} \ncarc{1}{0} \ncarc{2}{1}
\ncarc{1}{2} \ncarc{1}{3} \ncarc{3}{1} \ncarc{2}{3} \ncarc{3}{2}
\ncarc{4}{3} \ncarc{3}{4}\ncarc{3}{5} \ncarc{5}{3}\ncarc{5}{4}
\ncarc{4}{5}\ncarc{5}{6} \ncarc{6}{5}\ncarc{5}{7}
\ncarc{7}{5}\ncarc{6}{7} \ncarc{7}{6}\ncarc{7}{8}
\ncarc{8}{7}\ncarc{7}{9} \ncarc{9}{7} \ncarc{0}{3}\ncarc{3}{0}
\ncarc{2}{5}\ncarc{5}{2}\ncarc{4}{7}\ncarc{7}{4}\ncarc{6}{9}\ncarc{9}{6}
\psset{labelsep=-4.25ex}\nput{90}{0}{1}
\psset{labelsep=-4.25ex}\nput{90}{2}{3}
\psset{labelsep=-4.25ex}\nput{90}{4}{5}
\psset{labelsep=-4.25ex}\nput{90}{6}{7}
\psset{labelsep=-4.25ex}\nput{90}{1}{2}
\psset{labelsep=-4.25ex}\nput{90}{3}{4}
\psset{labelsep=-4.25ex}\nput{90}{5}{6}
\psset{labelsep=-4.25ex}\nput{90}{7}{8}
\end{psmatrix}
$$

\vspace{-1cm}

Finally the squared norms $\|Q_n\|_{\widetilde{W}}^2$ are $2\times2$ full matrices (not diagonal). Therefore, using Theorem \ref{Teo}, we get that an invariant distribution for this process is given by
\begin{equation}\label{ID2}
\mbox{\boldmath$\pi$}=
((\Pi_0\textbf{e}_2)^T;(\Pi_1\textbf{e}_2)^T;(\Pi_2\textbf{e}_2)^T;\cdots),
\end{equation}
where
\begin{equation}\label{Pn0N2}
\Pi_n=(\|Q_n\|_{\widetilde{W}}^2)^{-1},\quad n\geq0.
\end{equation}
We can use Theorem 8.1 of \cite{GdI2} to study the recurrence of the process, since our weight matrix has similar properties to that introduced in \cite{GdI2}. The Markov process that results from $P$ is never positive recurrent. If $-1<\beta\leq0$ then the process is null recurrent. If $\beta>0$ then the process is transient. For $-1<\beta\leq0$ the invariant distribution will be unique, as we mentioned at the end of Section \ref{IDs}.

The explicit expressions of $A_n$, $B_n$, $C_n$ and $\|Q_n\|_{\widetilde{W}}^2$ are too long to be displayed here. Instead, we will give all these data for the special case of $\alpha=\beta=0$ and $k=\frac{1}{2}$.  Similar expressions can be derived for other values of the parameters $\alpha, \beta$ and $k$. The coefficients of the three-term recurrence relation \eqref{ttrrst2} are, for $n\geq0$
\begin{align*}
A_n=&\left(\begin{array}{cc}
                                          \dfrac{(n+3)^2(2n^2+4n+1)(4n^2+18n+17)(16n^4+128n^3+348n^2+368n+117)}{2(n+2)(2n+3)(2n^2+8n+7)(4n^2+10n+3)(16n^4+160n^3+572n^2+860n+463)} & * \\
                                         \dfrac{2(n+3)^2(2n^2+4n+1)^2(4n^2+18n+17)}{(2n+3)(2n^2+8n+7)(4n^3+14n^2+9n+1)(16n^4+160n^3+572n^2+860n+463)}   & * \\
                                        \end{array}\right.\\
& \left.\begin{array}{cc}
                                          * & \dfrac{2(n+3)(2n^2+4n+1)(2n^2+12n+17)(4n^3+26n^2+49n+28)}{(n+2)(2n+3)(2n^2+8n+7)(4n^2+10n+3)(16n^4+160n^3+572n^2+860n+463)} \\
                                          * & \dfrac{(n+3)(2n^2+4n+1)(4n^3+26n^2+49n+28)(16n^4+128n^3+348n^2+368n+125)}{2(2n+3)(2n^2+8n+7)(4n^3+14n^2+9n+1)(16n^4+160n^3+572n^2+860n+463)} \\
                                        \end{array}\right),
\end{align*}
\begin{align*}
B_n=&\left(\begin{array}{cc}
                                         \dfrac{32n^8+384n^7+1928n^6+5256n^5+8450n^4+8148n^3+4577n^2+1365n+163}{(2n^2+4n+1)(2n^2+8n+7)(16n^4+96n^3+188n^2+132n+31)}  & * \\
                                         \dfrac{(n+2)(4n^2+10n+3)(32n^6+256n^5+760n^4+1040n^3+674n^2+186n+13)}{2(2n^2+4n+1)(2n^2+8n+7)(4n^3+14n^2+9n+1)(16n^4+96n^3+188n^2+132n+31)} & * \\
                                        \end{array}\right.\\
& \left.\begin{array}{cc}
                                          * & \dfrac{(4n^3+14n^2+9n+1)(32n^6+320n^5+1240n^4+2320n^3+2114n^2+834n+121)}{2(n+2)(2n^2+4n+1)(2n^2+8n+7)(4n^2+10n+3)(16n^4+96n^3+188n^2+132n+31)} \\
                                          * & \dfrac{(2n^2+6n+3)(4n^3+18n^2+21n+6)(4n^3+18n^2+21n+3)}{(2n^2+4n+1)(2n^2+8n+7)(16n^4+96n^3+188n^2+132n+31)} \\
                                        \end{array}\right),
\end{align*}
and for $n\geq1$
\begin{align*}
C_n=&\left(\begin{array}{cc}
                                         \dfrac{n(n+1)(4n^2+2n-3)(32n^6+192n^5+328n^4+32n^3-226n^2-4n+33)}{2(n+2)(2n+3)(2n^2-1)(4n^2+10n+3)(16n^4+32n^3-4n^2-20n+7)}  & * \\
                                         \dfrac{n(n+1)(4n^2+2n-3)(8n^4+32n^3+36n^2+8n-5)}{(2n+3)(2n^2-1)(4n^3+14n^2+9n+1)(16n^4+32n^3-4n^2-20n+7)} & * \\
                                        \end{array}\right.\\
 & \left.\begin{array}{cc}
                                          * &  \dfrac{n(4n^3+2n^2-7n+2)(8n^4+32n^3+36n^2+8n-5)}{(n+2)(2n+3)(2n^2-1)(4n^2+10n+3)(16n^4+32n^3-4n^2-20n+7)}\\
                                          * & \dfrac{n(4n^3+2n^2-7n+2)(32n^6+192n^5+392n^4+288n^3-34n^2-132n-43)}{2(2n+3)(2n^2-1)(4n^3+14n^2+9n+1)(16n^4+32n^3-4n^2-20n+7)} \\
                                        \end{array}\right).
\end{align*}
We see that all entries of the transition probability matrix \eqref{PP} are nonnegative rational functions in $n$, and it is very remarkable that the sum of the entries in each row of $P$ equals one.

The corresponding squared norms are the $2\times2$ non-diagonal matrices
\begin{align*}
\|Q_n\|_{\widetilde{W}}^2=&\left(\begin{array}{cc}
                                          \dfrac{(2n^2+4n+1)(16n^6+160n^5+628n^4+1212n^3+1173n^2+514n+79)}{(n+1)(2n+3)(4n^2+10n+3)^2(n+2)^3} & * \\
                                          -\dfrac{(2n+1)(2n+5)(2n^2+4n+1)}{2(n+1)(n+2)^2(4n^2+10n+3)(4n^3+14n^2+9n+1)}& * \\
                                        \end{array}\right.\\
 & \left.\begin{array}{cc}
                                          * & -\dfrac{(2n+1)(2n+5)(2n^2+4n+1)}{2(n+1)(n+2)^2(4n^2+10n+3)(4n^3+14n^2+9n+1)}\\
                                          * & \dfrac{(2n^2+4n+1)(16n^6+128n^5+388n^4+564n^3+417n^2+152n+22)}{(n+1)(n+2)(2n+3)(4n^3+14n^2+9n+1)^2} \\
                                        \end{array}\right).
\end{align*}

The unique (the process is null recurrent in this case) invariant distribution is given by \eqref{ID2} where $\Pi_n\textbf{e}_2$ can be calculated using \eqref{Pn0N2}
\begin{equation*}
\Pi_n\textbf{e}_2=\begin{pmatrix}
  \dfrac{2(n+1)^2(n+2)^2(2n+3)(4n^2+10n+3)(4n^2+14n+9)}{(2n^2+4n+1)(2n^2+8n+7)(16n^4+96n^3+188n^2+132n+31)} \\
 \dfrac{2(n+1)(n+2)(2n+3)(4n^3+14n^2+9n+1)(4n^3+22n^2+33n+8)}{(2n^2+4n+1)(2n^2+8n+7)(16n^4+96n^3+188n^2+132n+31)}\\
\end{pmatrix},\quad n\geq0.
\end{equation*}

If we denote by $\mbox{\boldmath$\pi$}^n_1$ and $\mbox{\boldmath$\pi$}^n_2$ the two components of $\Pi_n\textbf{e}_2$, we can study the behavior of the invariant distribution in Figure 4.1, a luxury we can afford since we have an analytic expression.

\bigskip

\bigskip

\begin{figure}[h]
\begin{center}\begin{minipage}{6cm}
\includegraphics[35,0][235,143]{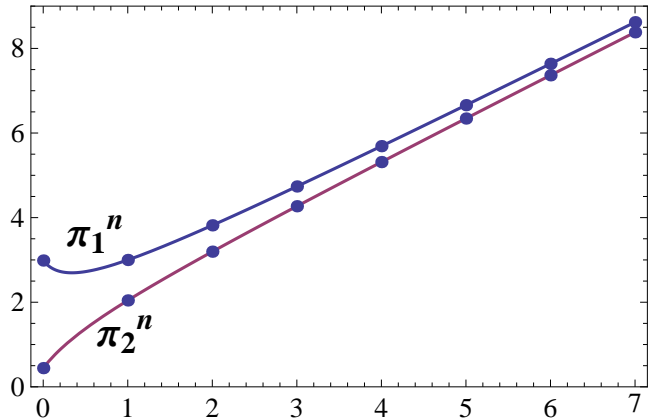}
\caption{$\alpha=0, \beta=0, k=0.5$}
\end{minipage}\end{center}
\end{figure}

\end{document}